\documentclass[twoside]{amsart}

\usepackage[hmargin={3cm,3cm},vmargin={3cm,3cm},includehead]{geometry}
\usepackage{amsmath,amssymb}
\usepackage{amsfonts}
\usepackage{mathrsfs}

\usepackage{dsfont}

\usepackage{epic}
\usepackage{eepic}

\newtheorem{theorem}{Theorem}
\newcommand{\ds}{\displaystyle}

\def\eqref#1{(\ref{#1})}

\newcommand{\bosn}{\textrm{\scriptsize $\boldsymbol{N}$}}

\newcommand{\bos}{\boldsymbol{a}}

\newcommand{\qalg}{\mathscr{A}}

\newcommand{\opr}{\stackrel{\textrm{def}}{=}}

\newcommand{\rop}{\mathtt{r}}
\newcommand{\e}{\boldsymbol{e}}

\begin{document}

\vspace{2cm}

\title[]{Tetrahedron equations and nilpotent subalgebras of $\mathscr{U}_q(sl_n)$.}%
\author{S. M. Sergeev}%
\address{Department of Theoretical Physics
\&
Department of Mathematics\\
Australian National University\\ Canberra ACT 0200\\ Australia}
\email{sergey.sergeev@anu.edu.au}



\begin{abstract}
A relation between $q$-oscillator $R$-matrix of the tetrahedron
equation and decompositions of Poinkar\'e-Birkhoff-Witt type bases
for nilpotent subalgebras
$\mathscr{U}_q(\mathfrak{n}_+)\subset\mathscr{U}_q(sl_n)$ is
observed.
\end{abstract}

\maketitle

The method of $q$-oscillator three-dimensional auxiliary problem,
derivation of three-dimensional $R$-matrix for the Fock space
representation of the algebra of observables, various
$q$-oscillator tetrahedron equations and details of $3D\to 2D$
dimension-rank transmutation may be found in \cite{ZTE}. In the
first part of this letter I remind the definition of
$q$-oscillator $R$-matrix of the tetrahedron equation.
Surprisingly, this $R$-matrix may be derived in a completely
different framework. In the second part of this letter it is shown
that matrix elements of $q$-oscillator $R$-matrix relate two
certain Poinkar\'e-Birkhoff-Witt type bases of the nilpotent
(positive roots) subalgebra $\mathscr{U}_q(\mathfrak{n}_+)\subset
\mathscr{U}_q(sl_3)$. The tetrahedron equation in this framework
follows from the completeness and irreducibility of
Poinkar\'e-Birkhoff-Witt type bases of the nilpotent subalgebra of
$\mathscr{U}_q(sl_4)$.

\noindent{\textbf{$Q$-oscillator $R$-matrix.}} Let elements
elements $\bos^+,\bos^-$ and $q^\bosn$ generate the $q$-oscillator
algebra $\qalg$:
\begin{equation}
\bos^+\bos^-=1-q^{2\bosn}\;,\quad
\bos^-\bos^+=1-q^{2\bosn+2}\;,\quad
q^{\bosn}\bos^{\pm}=\bos^{\pm}q^{\bosn\pm 1}\;.
\end{equation}
There is a special automorphism $x\to \rop_{123}^{} \cdot x\cdot
\rop_{123}^{-1}$ of $\qalg^{\otimes 3}$ defined by
\begin{equation}\label{mapping}
\left\{\begin{array}{l} \ds \rop_{123}^{}\cdot
q^{\bosn_2}\bos^\pm_1=\; (q^{\bosn_3}\bos^\pm_1 +
q^{\bosn_1}\bos^\pm_2\bos^\mp_3)\cdot
\rop_{123}^{},\\\\
\ds \rop_{123}^{} \cdot \bos^\pm_2  = (\bos^\pm_1\bos^\pm_3\,-
q^{1+\bosn_1+\bosn_3}\bos^\pm_2)\cdot \rop_{123}^{},\\\\
\ds \rop_{123}^{} \cdot q^{\bosn_2}\bos^\pm_3  =
(q^{\bosn_1}\bos^\pm_3
+q^{\bosn_3}\bos^\mp_1\bos^\pm_2)\cdot\rop_{123}^{}.
\end{array}\right.
\end{equation}
Here the indices of $q$-oscillator generators stand for the
components of the tensor product $\qalg^{\otimes 3}$. In the Fock
space (non-hermitian) representation
\begin{equation}\label{Fock}
\bos^-|0\rangle \;=\; 0\;,\quad |n\rangle \;=\;
\bos^{+n}|0\rangle\;,\quad q^{\bosn}|n\rangle \;=\; |n\rangle
q^n\;,\quad  \langle m|n\rangle = \delta_{n,m}
\end{equation}
the matrix elements of Eqs. (\ref{mapping}) provide a set of
recursion equations for the matrix elements of $\rop_{123}^{}$
with the unique solution
\begin{equation}\label{r-matrix}
\begin{array}{l}
\ds \langle m_1,m_2,m_3|\rop|n_1,n_2,n_3\rangle \;\opr\;
\rop_{m_1,m_2,m_3}^{\;n_1,\;n_2,\;n_3}\;=\\
\\
\ds \delta_{m_1+m_2,n_1+n_2} \delta_{m_2+m_3,n_2+n_3}
\frac{q^{(m_1-n_2)(m_3-n_2)}}{(q^2;q^2)_{m_2}}
P_{m_2}(q^{2n_1},q^{2n_2},q^{2n_3})\;,
\end{array}
\end{equation}
where $P_m(x,y,z)$ is defined recursively by
\begin{equation}
P_{m+1}(x,y,z)=(1-x)(1-z) P_m(\frac{x}{q^2},y,\frac{z}{q^2}) -
\frac{xz}{q^{2m}}(1-y) P_m(x,\frac{y}{q^2},z)\;,\quad
P_0(x,y,z)=1\;.
\end{equation}
Pohgammer's symbol is defined as usual,
$(q^2;q^2)_n=(1-q^2)(1-q^4)\cdots (1-q^{2n})$. Coefficients
$P_{m}(q^{2n_1},q^{2n_2},q^{2n_3})$ may be expressed in terms of
$q$-hypergeometric function $\phantom{}_2\varphi_1$. Matrix
$\rop_{123}$ is a square root of unity, $\rop_{123}^2=\mathds{1}$.
Note the symmetry property,
\begin{equation}
P_{m_2}(q^{2n_1},q^{2n_2},q^{2n_3})=\frac{(q^2;q^2)_{n_1}(q^2;q^2)_{n_3}}{(q^2;q^2)_{m_1}(q^2;q^2)_{m_3}}
P_{n_2}(q^{2m_1},q^{2m_2},q^{2m_3})\;,
\end{equation}
where $m_1+m_2=n_1+n_2$ and $m_2+m_3=n_2+n_3$.

Spectral parameters may be introduced as follows:
\begin{equation}
R_{123}\;=\;\left(-\frac{\lambda_1\mu_3}{q}\right)^{\bosn_2}\rop_{123}^{}
\left(\frac{\lambda_3}{\lambda_2}\right)^{\bosn_1}
\left(\frac{\mu_1}{\mu_2}\right)^{\bosn_3}\;,
\end{equation}
where $\mathbb{C}$-valued parameters $\lambda_j,\mu_j$ are
associated with $j$th component of a tensor power of $\qalg$.
Matrices $R_{ijk}$ satisfy the tetrahedron equation in
$\qalg^{\otimes 6}$,
\begin{equation}\label{TE}
R_{123}R_{145}R_{246}R_{356}=R_{356}R_{246}R_{145}R_{123}\;.
\end{equation}
Due to the $\delta$-symbols structure of (\ref{r-matrix}), the
spectral parameters may be removed from the tetrahedron equation,
Eq. (\ref{TE}) is equivalent to the constant tetrahedron equation.
However, the constant $\rop$-matrix provides spectral parameter
dependent solutions for the Yang-Baxter equation
\begin{equation}
R_{s_1,s_2}(u)R_{s_1,s_3}(uv)R_{s_2,s_3}(v)=R_{s_2,s_3}(v)R_{s_1,s_3}(uv)R_{s_1,s_2}(u)\;.
\end{equation}
For instance, in the lowest rank case of affine
$\mathscr{U}_q(\widehat{sl}_2)$, matrix elements of
two-dimensional  $R$-matrix of the Yang-Baxter equation in spin
$s_1$ $\times$ spin $s_2$ evaluation representations are given by
\begin{equation}\label{6v}
\langle m_1,m_2|R_{s_1,s_2}(u)|n_1,n_2\rangle
\;=\;\sum_{n_3,m_3=0}^{\infty} u^{m_3}
\rop_{s_1+m_1,s_2+m_2,m_3}^{\;s_1+n_1,\;s_2+n_2,\;n_3}
\rop_{s_1-m_1,s_2-m_2,\;n_3}^{\;s_1-n_1,\;s_2-n_2,m_3}\;,
\end{equation}
where $m_j,n_j\in \{-s_j,-s_j+1,\cdots, s_j-1,s_j\}$, and the sum
converges if $|q|<1$ and $|u|<1$. In particular, formula
(\ref{6v}) gives the six-vertex $R$-matrix for $s_1=s_2=1/2$.

\vspace{1cm}

\noindent{\textbf{Nilpotent subalgebras.}} Let
$\mathscr{B}_n=\mathscr{U}_q(\mathfrak{n}_+)$ be the nilpotent
subalgebra of
$\mathscr{U}_q(sl_{n+1})=\mathscr{U}_q(\mathfrak{n}_-)\otimes
\mathscr{U}_q(\mathfrak{h})\otimes \mathscr{U}_q(\mathfrak{n}_+)$.
It is generated by elements $\e_1,\e_2,\dots,\e_n$ satisfying
$q$-Serre relations:
\begin{equation}\label{borel}
\e_i^2\e_{i\pm 1}+\e_{i\pm 1}\e_i^2 \;=\; (q^{-1}+q) \e_i\e_{i\pm
1} \e_i\quad\textrm{and}\quad \e_i\e_j=\e_j\e_i\;\;\;
\textrm{if}\;\;\; |i-j|>1\;.
\end{equation}
Consider for a moment $\mathscr{B}_2$ and define
\begin{equation}
\e_{12}^{}=\frac{\e_1\e_2-q\e_2\e_1}{q^{-1}-q}\quad
\textrm{and}\quad
\widetilde{\e}_{12}=\frac{\e_2\e_1-q\e_1\e_2}{q^{-1}-q}\;.
\end{equation}
It is well known \cite{Lusztig}, there are many ways to define a
Poinar\'e-Birkhoff-Witt type basis in $\mathscr{B}_n$. For
instance, two following bases of $\mathscr{B}_2$,
\begin{equation}
\biggl\{\e_2^{n_3}\e_{12}^{n_2}\e_1^{n_1},\;\;n_j\geq 0\biggr\}
\quad\textrm{and}\quad
\biggl\{\e_1^{m_1}\widetilde{\e}_{12}^{m_2}\e_2^{m_3},\;\;m_j\geq
0\biggr\}\;,
\end{equation}
are complete, irreducible and therefore equivalent. The
equivalence of the bases is given by the decomposition
\begin{equation}\label{decomp}
\frac{\e_1^{m_1}}{[m_1]!} \frac{\widetilde{\e}_{12}^{m_2}}{[m_2]!}
\frac{\e_2^{m_3}}{[m_3]!} \;=\; \sum_{\{n\}}\;
\rop_{m_1,m_2,m_3}^{\;n_1,\;n_2,\;n_3} \;
\frac{\e_2^{n_3}}{[n_3]!} \frac{\e_{12}^{n_2}}{[n_2]!}
\frac{\e_1^{n_1}}{[n_1]!}
\end{equation}
where we use for shortness
\begin{equation}
[n]=q^{-n}-q^n\;,\quad [n]!\;=\;[1][2]\cdots [n]\;=\;
q^{-n(n+1)/2}(q^2,q^2)_n\;.
\end{equation}
\begin{theorem}
The coefficients $\rop_{m_1,m_2,m_3}^{\;n_1,\;n_2,\;n_3}$ in
(\ref{decomp}) are given exactly by the formula (\ref{r-matrix}).
\end{theorem}
\noindent Sketch proof is the following. Coefficients
$\rop_{m_1,m_2,m_3}^{\;n_1,\;n_2,\;n_3}$ in (\ref{decomp}) satisfy
a set of recursion relations following from the structure of
bases. For instance, a simple identity following from
(\ref{borel})
\begin{equation}
\e_1^{m_1+1}\widetilde{\e}_{12}^{m_2} \e_2^{m_3}\;=\;
q^{m_2}[m_3]\e_1^{m_1}\widetilde{\e}_{12}^{m_2}\e_2^{m_3-1}\e_{12}^{}
+q^{m_2+m_3}\e_1^{m_1}\widetilde{\e}_{12}^{m_2}\e_2^{m_3}\e_1^{}
\end{equation}
provides
\begin{equation}
(1-q^{2m_1+2})\rop_{m_1+1,m_2,m_3}^{\;\;\;n_1,\;n_2,\;n_3}\;=\;
\rop_{m_1,m_2,m_3-1}^{\;n_1,\;n_2-1,\;n_3} (1-q^{2n_2}) + q^{m_3}
\rop_{\;\;m_1,m_2,m_3}^{n_1-1,\;n_2,\;n_3}q^{n_2}(1-q^{2n_1})\;.
\end{equation}
In terms of $q$-oscillators in the basis (\ref{Fock}) it
corresponds to
\begin{equation}
\bos_2^-\rop_{123}^{} =
\bos_3^+\rop_{123}^{}\bos_2^-+q^{\bosn_3}\rop_{123}^{}
q^{\bosn_2}\bos_1^-\;,
\end{equation}
what is simple consequence of (\ref{mapping}). The complete
collection of all such recursion relations for (\ref{decomp}) is
equivalent to the system (\ref{mapping}).\hfill$\blacksquare$

Consider next $\mathscr{B}_3$ and equivalence of the complete
irreducible bases
\begin{equation}
\biggl\{ \frac{\e_3^{n_6} \e_{23}^{n_5} \e_{123}^{n_4} \e_2^{n_3}
\e_{12}^{n_2} \e_1^{n_1}}{[n_6]![n_5]!\cdots [n_1]!}\;,\;\;n_j\geq
0\biggr\}\;,
\end{equation}
where
\begin{equation}
\e_{12}=\frac{\e_1\e_2-q\e_2\e_1}{[1]}\;,\quad \e_{23}=
\frac{\e_2\e_3-q\e_3\e_2}{[1]}\;,\quad
\e_{123}=\frac{\e_1\e_{23}-q\e_{23}\e_1}{[1]}\;,
\end{equation}
and
\begin{equation}
\biggl\{ \frac{\e_1^{m_1} \widetilde{\e}_{12}^{m_2}
\widetilde{\e}_{123}^{m_4} \e_2^{m_3} \widetilde{\e}_{23}^{m_5}
\e_3^{m_6}}{[m_1]![m_2]!\cdots [m_6]!}\;,\;\;m_j\geq 0\biggr\}\;,
\end{equation}
where
\begin{equation}
\widetilde{\e}_{12}= \frac{\e_2\e_1-q\e_1\e_2}{[1]}\;,\quad
\widetilde{\e}_{23}= \frac{\e_3\e_2-q\e_2\e_3}{[1]}\;,\quad
\widetilde{\e}_{123} =
\frac{\widetilde{\e}_{23}\e_{1}-q\e_{1}\widetilde{\e}_{23}}{[1]}\;.
\end{equation}
The decomposition matrix $T_{m_1,\dots,m_6}^{\;n_1,\dots,\;n_6}$
for these bases may be constructed in two alternative ways,
\begin{equation}
T_{12\dots
6}^{(1)}=\rop_{123}\rop_{145}\rop_{246}\rop_{356}\quad\textrm{and}\quad
T_{12\dots 6}^{(2)}=\rop_{356}\rop_{246}\rop_{145}\rop_{123}\;.
\end{equation}
Due to the uniqueness of the decomposition, $T_{12\dots
6}^{(1)}=T_{12\dots 6}^{(2)}$, what is the tetrahedron equation.

\vspace{1cm}

\noindent{\textbf{``Coherent'' states.}} The result of equivalence
of (\ref{mapping}) and (\ref{decomp}) may be written in a
remarkable simple form in terms of ``coherent states''. Let
\begin{equation}
\psi(u)\;=\;\sum_{n=0}^\infty \frac{u^n}{[n]!}\;=\;
\sum_{n=0}^\infty q^{n(n+1)/2} \frac{u^n}{(q^2;q^2)_n}\;,\quad
\psi(u/q)-\psi(qu)=u\psi(u)\;.
\end{equation}
Then (\ref{mapping}) and (\ref{decomp}) may be combined together
into the basis-invariant form in $\qalg^{\otimes 3}\otimes
\mathscr{B}_2$:
\begin{equation}\label{psi-vectors}
\rop_{123}^{}\cdot \psi(\e_2^{}\bos_3^+) \psi(\e_{12}^{}\bos_2^+)
\psi(\e_1^{}\bos_1^+) |0\rangle \;=\; \psi(\e_1^{}\bos_1^+)
\psi(\widetilde{\e}_{12}^{}\bos_2^+) \psi(\e_2^{}\bos_3^+)
|0\rangle\;,
\end{equation}
where $|0\rangle$ is the total Fock vacuum for $\qalg^{\otimes
3}$, $\rop_{123}^{}$ is the $q$-oscillator $R$-matrix of
(\ref{mapping}), and the auxiliary coefficients $\e_1,\e_2$ $\in$
$\mathscr{B}_2$. Presumably, relation (\ref{psi-vectors}) is the
master equation for three-dimensional variant of the vertex-IRF
duality.


\end{document}